\documentclass[twoside,12pt,leqno]{article}
\usepackage[margin=20mm,bmargin=35mm,tmargin=30mm]{geometry}
\usepackage[english]{babel}
\usepackage[utf8]{inputenc}
\usepackage{graphicx}
\usepackage{mathtools}
\usepackage{epstopdf}
\usepackage{pifont}
\usepackage{bbding}
\usepackage{units}
\usepackage{amsthm}
\usepackage{amsfonts}
\usepackage{titlesec}
\usepackage{fancyhdr}
\theoremstyle{definition}
\newtheorem{definition}{Definition}

\title{\vspace{-20mm}\Large Four redundant axiomatics}
\author{\large Anton Cedilnik\\ \normalsize University of Ljubljana}
\date{}

\makeatletter
\renewcommand{\thesection}{\@arabic\c@section.\hspace*{-2mm}}
\renewcommand{\section}{\@startsection{section}{1}{\z@}%
                       {-3.5ex \@plus -1ex \@minus -.2ex}%
                       {2.3ex \@plus.2ex}%
                       {\normalfont\large\centering\bfseries}}
\makeatother

\begin{document}

\maketitle
\thispagestyle{empty}

\pagestyle{fancy}

\begin{abstract}
  \footnotesize
  In the paper, we shorten four axiom systems which are redundant in the mathematical literature.

  \vspace{5mm}
  

  \vspace{5mm}

  \noindent \textit{AMS math. subj. class. 08-01}
\end{abstract}

\baselineskip=14.5pt

\smallskip

\section{}
Axiomatics has been an essential method in mathematics from Euclides' times up to today. The fewer axioms and the more consequences, the more successful is the theory. Of course, this is not the only criterion of successfulness, much more important is the applicability of the theory, but this will not be of our interest here. Since there should be as few assumptions as possible, we wish no assertion which could be derived from other axioms to be inserted into the axiom system, i.e., that the system is not redundant. Although the purpose of the paper is mainly didactic, we believe that the fourth example also has some scientific value.

Out of a desire not to waste time on deriving almost trivial consequences, we often consciously or at least without careful consideration set up a system of axioms which can sometimes be reduced with rather small corrections.

The first example of such system will be the definition of metric. The following is from the standard sources [1], p. 28, and [5], p. 9 (with unimportant stylistic changes):

\begin{definition}
\textbf{Metric space} \ $(\mathcal{M},d)$ \ is an ordered pair where 
\begin{itemize}
\item[\ding{226}] $\mathcal{M}$ \, is set,
\item[\ding{226}] $d$ \, is a \textbf{metric} on $\mathcal{M}$\,, \, i.e., a function $d : \mathcal{M}^2 \rightarrow \mathbb{R}$,\\ such that for any \, $(x , y , z) \, \in \, \mathcal{M}^3 \,,$\, the following holds:
\end{itemize}
\noindent \fbox{\begin{minipage}{41em}
\vspace{0.1cm}
\textbf{M1.} \ \ $d(x , y) \, \geq \, 0$   \\
\textbf{M2.} \ \ $d(x , y) \, = \,  0$ \  \ $\Leftrightarrow$ \  \ $x = y$  \\     
\textbf{M3.} \ \ $d(x , y) \, = \, d(y , x) $   \\ 
\textbf{M4.} \ \ $d(x , y) \, \leq \, d(x , z) + d(z , y) $   
\end{minipage}}
\vskip 4mm
But in fact two axioms are enough:
\\\\
\noindent \fbox{\begin{minipage}{41em}
\vspace{0.1cm}
\textbf{M2.} \ \ \ \ \ \ \  $d(x , y) \, = \,  0$ \  \ $\Leftrightarrow$ \  \ $x = y$  \\  
\textbf{M4}$^{\rm{var}}.$ \ \ \ \ $d(x , y) \, \leq \, d(x , z) + d(y , z) $   
\end{minipage}}
\vskip 4mm
\noindent We prove M1 and M3 by inserting firstly \ $y = x$ \ and secondly \ $z = x$ \ into \ M4$^{\rm{var}}.$ \ M4 \\ then follows from \ M4$^{\rm{var}}$ and M3.

\end{definition}

A similar example is the definition of the norm. Let us quote [1], p. 92, and [5], p. 95:\\
\begin{definition}
A \textbf{normed vector space} \  $(\,\mathcal{N}, ||.||\,)$ \  is an ordered pair where

\begin{itemize}
\item[\ding{226}]  $\,\mathcal{N}$ \ is a vector space over \ $\mathbb{F} \, \in \{ \mathbb{R}, \mathbb{C} \}$\,,   
\item[\ding{226}]   $||.||$ \ a \textbf{norm} on \, $\mathcal{N}$\,, \, i.e., \, a function \ $||.||: \ \mathcal{N} \  \rightarrow  \ \mathbb{R}\,,\ ||.||\,:\, x \  \mapsto \  ||x||\,,$   \\
such that for any \ $(\lambda , x , y) \, \in \, \mathbb{F} \times \mathcal{N}^2$ \ the following holds:  \ 
\end{itemize}
\end{definition}
\noindent \fbox{\begin{minipage}{41em}
\vspace{0.1cm}
\textbf{N1.} \ \  $ ||x|| \, \geq \, 0 $ \\  
\textbf{N2.} \ \  $ ||x|| \, = 0  \ \ \Leftrightarrow \ \ x = 0 $ \\ \textbf{N3.} \ \  $ ||\lambda x|| \, = \, |\lambda| \, ||x|| $ \\  
 \textbf{N4.} \ \  $ || x + y|| \, \leq \, ||x|| + ||y|| $ 
\end{minipage}}
\vskip 4mm
Here too, the axiom system can be considerably reduced:\\\\
\noindent \fbox{\begin{minipage}{41em}
\vspace{0.1cm}
\textbf{N2}$^{\rm{var}}.$\ \ \ \ \ \ \ $ ||x|| \, = \, 0 \ \Rightarrow \ x = 0 $ \\  
\textbf{N3}$^{\rm{var}}.$\ \ \ \ \  \ \  $ ||\lambda x|| \, \leq \, |\lambda| \, ||x|| $ \\  
\textbf{N4.} \ \ \ \ \ \ \ \ \  $ || x + y|| \, \leq \, ||x|| + ||y|| $ 
\end{minipage}}
\vskip 4mm
\noindent In addition to the formal, the substantive advantage of the axiom N3$^{\rm{var}}$ before N3 is that it is semantically consistent with N4 and also with the submultiplicativity of the norm, if we upgrade the normed vector space into a normed algebra: \ $||x \cdot y ||\, \leq \, ||x||\, ||y|| \,.$\\

From N3$^{\rm{var}}$ follows for \ $\lambda = 0$ \, : \ $||0|| \leq 0$\,; \ and from N4 for \ $x = y = 0\,,$\ the inverse relation follows. So \ $||0|| = 0$ \ and N2 holds. Then we get N1 by derivation

$$ 0 = \frac{1}{2} ||0|| = \frac{1}{2} ( || x + (-1) x || \, \leq  \, \frac{1}{2} (||x ||) + |-1| \, ||x|| ) = ||x|| \,. $$
Another derivation (for $\lambda \neq 0$\,)
$$ |\lambda| \, ||x|| = |\lambda| \left\| \frac{1}{\lambda} \cdot \lambda x  \right\| \, \leq \, |\lambda| \left| \frac{1}{\lambda}  \right| ||\lambda x|| \, = \, ||\lambda x|| \,  \leq \, |\lambda|\, ||x||$$
confirms  N3.
\section{}
\ \ \ \ \ In the basic courses of higher mathematics in technical as well as in some social science fields of study, the concept of vector space is included in the program. The goals are usually the three-dimensional vector space of geometric vectors and the Euclidean spaces $\mathbb{R}^n$, sometimes also the spaces $\ell_2$ and $\mathcal{L}_2$ . As a rule, therefore, the program contains real vector spaces with inner product, but not the abstract vector spaces. The system of axioms is therefore usually as follows ([2], p.~309, 310, 312).

\begin{definition}
\textbf{Inner product space} \ $(\mathcal{S}, +, \cdot,\langle.\,,. \rangle )$ \  is an ordered quadruple where
\begin{itemize}
\item[\ding{226}] $\mathcal{S} \neq$ \O \   is a set of elements called \textbf{vectors},
\item[\ding{226}] $+ : \mathcal{S}^2 \rightarrow \, \mathcal{S}$ \  is a binary operation called \textbf{addition},
\item[\ding{226}] $\cdot : \mathbb{R} \times \mathcal{S} \, \rightarrow \mathcal{S} $ \ is a binary operation called \textbf{multiplication by number},
\item[\ding{226}]  $\langle . \, , . \rangle \, : \, \mathcal{S}^2 \rightarrow \, \mathbb{R} \,$ \ is a binary operation called \textbf{inner product}. 
\end{itemize}
The following axioms must hold for any \ $(\lambda, \mu , x , y , z) \, \in \, \mathbb{R}^2 \times \mathcal{S}^3$:\\\\
\noindent \fbox{\begin{minipage}{41em}
\vspace{0.1cm}
\textbf{S1.} \ \  $ x + y = y 
+ x$ \\  
\textbf{S2.} \ \  $ (x + y) + z = x + (y + z) $ \\ 
\textbf{S3.} \ \   $\exists \, 0 \, \in \, \mathcal{S} \ : \ x + 0 \ = \ x  $ \\
\textbf{S4.} \ \  $ \exists \, (-x)\, \in \, \mathcal{S} \ : \ x + (-x) \ = \ 0    $ \\
\textbf{S5.} \ \  $ (\lambda + \mu) \cdot x \ = \ \lambda \cdot x + \mu \cdot x  $ \\
\textbf{S6.} \ \  $ \lambda \cdot (x + y) \ = \ \lambda \cdot x + \lambda \cdot y  $ \\
\textbf{S7.} \ \  $ \lambda \cdot (\mu \cdot x) \ = \ (\lambda \mu) \cdot x $ \\
\textbf{S8.} \ \  $ 1 \cdot x  \ = \ x$ \\
\textbf{S9.} \ \  $ \langle x , y \rangle =  \langle y , x \rangle $ \\
\textbf{S10.} \   $ \langle x + y,z \rangle =  \langle  x , z \rangle + \langle  y , z \rangle   $ \\
\textbf{S11.} \   $  \langle \lambda  \cdot x, y \rangle =    \lambda \langle  x , y \rangle   $ \\
\textbf{S12.} \   $ x \neq 0 \ \ \Rightarrow \ \     \langle x , x \rangle > 0$ 
\end{minipage}}
\vskip 4mm
From the didactic point of view, this system of axioms is complete nonsense. A student who, of course, does not know the background of these axioms, has to learn them by heart, which requires a great deal of effort, especially since the order of axioms is to some extent important. Since various trivialities are derived from them, such as the uniqueness of element 0 and opposite element and equality of type $0 \cdot x = 0$\,, \ the student inadvertently wonders why not take all these facts as axioms at all; if there are already 12 initial ones, what if there are a few more! Even the information that the first eight axioms define a vector space (and the first four an Abelian group) will mean nothing to the student.

The way out of this dilemma is in the realization that this system of axioms is severely redundant. Only four axioms are enough!
\\\\
\noindent \fbox{\begin{minipage}{41em}
\vspace{0.1cm}
\textbf{S9.} \ \  $ \langle x , y \rangle =  \langle y , x \rangle $ \\
\textbf{S10.} \   $ \langle x + y,z \rangle =  \langle  x , z \rangle + \langle  y , z \rangle   $ \\
\textbf{S11.} \   $  \langle \lambda  \cdot x, y \rangle =    \lambda \langle  x , y \rangle   $ \\
\textbf{S13.} \   $ \langle x , x \rangle +  \langle y , y \rangle = 2 \langle x , y \rangle \ \ \Rightarrow \ \ x = y  $
\end{minipage}}
\end{definition}
Only S13 is a completely new axiom, but it has a fairly obvious meaning which is easily understood by students:
$$ a^2 + b^2 = 2ab \ \ \Rightarrow \ \ a^2 - 2ab + b^2 \, = \, (a-b)^2 \, = \, 0 \ \  \Rightarrow \ \ a = b   $$
The trivial finding is that from the original system of axioms follows the shrunken system. Let us prove that this also holds in the opposite direction!

From S9 and S10 follows the identity
$$ \langle w + x , y + z \rangle = \langle w , y \rangle + \langle x , y \rangle + \langle w , z \rangle + \langle x , z \rangle \, ,$$
and from S9 and S11 also
$$ \langle \lambda \cdot x , \mu \cdot y \rangle \,  = \, \lambda \mu \langle x , y \rangle  \, .$$

The identities S1, S2, S5, S6, S7 and S8 are all provable by the same method. If we denote with $L$ the left side and with $R$ the right side of any of these identities, we calculate $\langle L , L \rangle $ ,  $\langle R , R \rangle $ and $\langle L , R \rangle $ \, and use the axiom S13. By the same procedure we prove \, $0 \cdot x = 0 \cdot y $ \, for any elements $x , y$\,.

Let's fix some element $w$ and calculate for any $x$:
$$x + 0 \cdot w  = 1 \cdot x + 0 \cdot x  = (1 + 0) \cdot x = x \, , $$\\
and S3 is proven. We prove S4 in a similar way:
$$ x + (-1) \cdot x = 1 \cdot x + (-1) \cdot x = (1 + (-1)) \cdot x = 0 \cdot x = 0 \, . $$
At this point, we can already define the subtraction.

Only S12 remained. First, let's calculate this for any element $x$ :
$$ \langle 0 , x \rangle = \langle 0 \cdot  x, x \rangle = 0 \langle x , x \rangle = 0 \,. $$
If \ $\langle x , x \rangle = 0$ \, , \ due to \ $\langle x , x \rangle + \langle 0 , 0 \rangle = 2 \langle x , 0 \rangle  $ \ we conclude from S13 that \ $x = 0$\,. \ 
So: \, $x \neq 0 \ \Rightarrow \ \langle x , x  \rangle \neq 0 $\,. \ Now suppose that for selected \ $y$ \, and \, $z$
$$ \langle y , y\rangle > 0 >  \langle z , z \rangle \,.$$
Then \ $ \langle \lambda \cdot y - z , \lambda \cdot y - z \rangle \ = \ 0 \,, $ \  if
$$\lambda = \left( \langle y , z \rangle + \sqrt{\langle y , z  \rangle^2 + \langle y , y \rangle | \langle z , z \rangle | }    \right)  \langle y , y \rangle^{-1} \,. $$
Hence\, $z = \lambda \cdot y  $\ and therefore \ $\langle z , z \rangle = \lambda^2 \langle y , y \rangle > 0  \,,$ \ which is a contradiction.

Thus, we have found out that the expressions $ \langle y , y \rangle$ for $y \neq 0$ are either all positive or all negative. In case they are negative, we replace the original inner multiplication $x , y \mapsto  \langle x , y \rangle $ with the new $x , y \mapsto - \langle x , y \rangle $ because all axioms still apply to it.\\

Of course, this means that the original S12 does not actually follow from the reduced axiom system. It is necessary to add to the shrunken system the agreement that the scalar product is that of $x , y \mapsto \pm \langle x , y \rangle $ , which is positive definite. Since both variants are already guaranteed by the other assumptions, this agreement need not be taken as a new axiom.

Three brief comments on the new set of axioms. The first point of interest is that all information about vectors is given exclusively through the inner product. Kind of like this: vectors are ``mysterious" \ objects  that can only be known through inner multiplication.

Similar work that we are doing here can be found in the paper [4] about the fact that we need six axioms for the definition of a vector space (Definition 3 has defined the vector space by S1 - S8). The paradox is that for a much richer structure, i.e. for the inner product space, only our four axioms are sufficient.

It is also important to note that these axioms could be met by scalars from any ring, since the order does not directly play any role, which gives ample opportunity to a researcher to generalize the concept of inner product space (as it is for example Hilbert module). Axiom S9 should be corrected somehow in this style:\\
\textbf{S9$^{\rm{var}}$}. \ \  $ \langle w , x \rangle =  \langle y , z \rangle  \  \Rightarrow  \ \langle x , w \rangle =  \langle z , y \rangle   \,,$ \\
and the axiom S12 should also look something like this:\\
\textbf{S12$^{\rm{var}}$}. \ \  $ \forall x , y \ \exists \, \lambda \,  :  \, \langle y , x \rangle  \ =  \lambda \langle x , x \rangle   \,.$ 
\section{}
\ \ \ \ \ When we talk about algebra with scalar involution, we usually assume the definition in [3]. Here we will write it in the same style as the previous three.
\begin{definition}
\textbf{Algebra with scalar involution } \ $(\mathcal{C}, \cdot , e , \ast)$ \ is an ordered quadruple where
\begin{itemize}
\item[\ding{226}] $\mathcal{C} \neq \{0\}$ \  is a vector space over a field \,$\mathbb{F}\,$,
\item[\ding{226}]   $(\mathcal{C}, \cdot , e)$  \ is a (possibly non-associative) algebra with a multiplication $\cdot$ and unit $e$,
\item[\ding{226}]  $\ast \, : \mathcal{C} \mapsto \mathcal{C} \, ,\ \ast : x \mapsto x^\ast $\, , \ is a map for which the following relations hold for any \\ \ $(\lambda, x , y) \ \in \mathbb{F} \times \mathcal{C}^2$\,:
\end{itemize}
\noindent \fbox{\begin{minipage}{41em}
\vspace{0.1cm}
\textbf{C1.} \ \  $ (x + y)^\ast = x^\ast + y^\ast $ \\  
\textbf{C2.} \ \  $ (\lambda x )^\ast  = \lambda x^\ast  $ \\ 
\textbf{C3.} \ \  $ (x \cdot y)^\ast = y^\ast \cdot x^\ast $ \\   
\textbf{C4.} \ \  $ x^{\ast \ast}  = x$ \\
\textbf{C5.} \ \  $ x + x^\ast \, \in \, \mathbb{F}\{e\} $ \\
\textbf{C6.} \ \  $ x \cdot x^\ast \, \in \, \mathbb{F}\{e\} $ 
\end{minipage}}
\vskip 4mm
 The definition of such an algebra necessarily includes the definition of a functional
$$ t(x) e \, : = \, - (x +
 x^\ast)\,.  $$
 
 We will replace all six axioms with only two.\\\\
 \noindent \fbox{\begin{minipage}{41em}
\vspace{0.1cm}
\textbf{C3.} \ \  $ (x \cdot y)^\ast = y^\ast \cdot x^\ast $ \\  
\textbf{C5.} \ \  $ x + x^\ast \, \in \, \mathbb{F}\{e\} $
\end{minipage}}
\vskip 4mm

 It immediately turns out, however, that these two axioms are not enough!\\
 
 Let us look at a one-dimensional algebra \ $(\mathbb{F}\{ e\}, e \cdot e = e )$ \ with mapping \ $ x \mapsto x^\ast $ \, ,  \ which defines a function \ $s : \mathbb{F} \rightarrow \mathbb{F}  $ according to the equation 
 $$ \forall \lambda \, \in \, \mathbb{F} \, : \, s(\lambda) e : = (\lambda e)^\ast  . $$
 Axiom C3 (and C5 trivially) holds exactly when the function $s$ has the property
  $$ \forall (\lambda,\mu) \, \in \, \mathbb{F}^2 \, : \, s(\lambda \mu)  = s(\lambda) \, s(\mu) . $$
\end{definition}
There are 3 possibilities:
\begin{itemize}
    \item $s \equiv  0 $ \  or \ $s \equiv 1 \,,$
    \item $s(0) = 0 $\, , \ $s(1) = 1$ \ and the restriction of $s$ is endomorphism of the group \ $(\mathbb{F}^{\times},\cdot)$\,.
\end{itemize}
If axioms C2 and C4 are to hold, the function $s$ should be the identity, which is obviously not necessary!

It is quite clear that in this case the submitted shrunken system should be replaced by only one axiom:\\
{\bf C$2^{\rm{var}}.$} \ \ $(\lambda e)^\ast = \lambda e $\,.\\

Notwithstanding this counterexample, let us insist with only two axioms! The algebra $(\mathcal{C}, \cdot , e , \ast)$ should therefore satisfy the axioms C3 and C5, together with a necessary demand $\rm{dim} \mathcal{C} \geq 2$\,.\\

Subtract the equations \ $x + x^\ast = -t(x) e$ \, and \, $x^\ast + x^{\ast \ast} = -t(x^\ast) e \,.$\ Then :
\begin{equation}
 \forall \, x \, \in \, \mathcal{C} \, : \ x^{\ast \ast} = x + [t(x) - t(x^\ast)] e \, .  
\end{equation}
Next conclusion: since \ $\lambda e + (\lambda e)^\ast = -t(\lambda e)e \ $ then 
\begin{equation}
\forall \, \lambda \, \in \,  \mathbb{F} : \ (\lambda e)^\ast \, = \,- [\lambda + t(\lambda e)]e \, \in \, \mathbb{F}\{e\}\,,
\end{equation}
\begin{equation}\forall \, (\lambda,x) \, \in \, \mathbb{F} \times \mathcal{C} : \ (\lambda x)^\ast  = (\lambda e \cdot x)^\ast = x^\ast \cdot (\lambda e)^\ast = -[\lambda + t(\lambda e)] x^\ast\,. \end{equation}
Because of $x^\ast + x^{\ast \ast} = -t(x^\ast)e$ \ and (1) and (2) it holds:

$x^\ast \, \in \, \mathbb{F}\{e\} \ \Rightarrow \  x^{\ast \ast} \, \in \mathbb{F}\{e\} \ \Rightarrow \  x \in \mathbb{F}\{e\}  \  \Rightarrow \  x^\ast \, \in \, \mathbb{F}\{e\}\,.   $\\
Therefore \ $x^\ast \, \in \, \mathbb{F}\{e\} \ $  if and only if \ $x \, \in \mathbb{F}\{e\}$\,.

We already know from (2) that $e^\ast \, \in \,  \mathbb{F}\{e\}\,;\, $ say, $ \, e^\ast = \varepsilon e \ . $ \ Then: 
$$ \varepsilon e = e^\ast = (e \cdot e)^\ast = e^\ast \cdot  e^\ast  =  \varepsilon^2 e \ , $$
$$(\varepsilon - \varepsilon^2)\,e = 0 \ ,
$$
which means that either \ $e^\ast = 0$ \ or  \ $e^\ast = e$\,. \ But because of C5, for any  \ $x \in \mathcal{C} \backslash \mathbb{F} \{e\} $ \ holds a similar equality: 
$$0 \neq x^\ast \ = \ (e \cdot x)^\ast \ = \ x^\ast \cdot e^\ast ,$$
which is possible only if $e^\ast = e$\ .

From the following calculation, using (3),
$$ -t(\lambda x)e - \lambda x \  = \ (\lambda x)^\ast \ = \  x^\ast \cdot ( \lambda e)^\ast$$
$$=[-t(x)e - x] \cdot [-t(\lambda e)e - \lambda e] \ = \ t(x) [t (\lambda e) + \lambda]e + [t(\lambda e)+ \lambda] x  $$
by combining the beginning and the end and assuming linear independence of the pair \ $(e , x)$ \ , \ we get : \ $t(\lambda e) = -2 \lambda$\,. \ Then from (3) we find:
$$\forall \,(\lambda, x) \, \in \mathbb{F} \times  \mathcal{C} \, : \ (\lambda x)^\ast \, = \, \lambda x^\ast \ ,$$ 
which is C2. Consequences: C7 and $t(e) = -2$ \ and
$$\forall \,(\lambda, x) \, \in \mathbb{F} \times \mathcal{C} \, : \ t(\lambda x) \, = \, \lambda t (x)\  . $$\\

Let's derive further!
$$ x \cdot x^\ast \  = \ x \cdot (x + x^\ast - x ) \ = \ x \cdot (-t(x)e - x) \ = \ -t(x)x - x \cdot x \ .    $$
So, the following identity holds always:
\begin{equation}
  \forall \, x \ \in \, \mathcal{C}\, : \, x \cdot x + t(x)x + x \cdot x^\ast \ = \ 0 \,.    
\end{equation}

If from the equation \ $(x + y) + (x + y)^\ast \, = \, -t(x + y)e \,  \  $  we subtract the equations \ $x + x^\ast \, = \, -t(x)e \  $ and \ $y + y^\ast \, = \, -t(y)e $ \ , \ then :
\begin{equation} \ \ \ \  (x + y)^\ast = [t (x) + t(y) - t(x + y)]e + x^\ast + y^\ast \,. \end{equation} 

\noindent Use (5) together with such an element $z$, \  such that the pair \ $(e , z)$  \ is linearly independent:
\vskip 2mm
\noindent $ [t(x) + t (y) - t (x + y)]z^\ast + z^\ast \cdot \  x^\ast + z^\ast \cdot \  y^\ast = z^\ast \cdot \  (x + y)^\ast = [(x + y) \cdot \  z]^\ast  = \\ (x \ \cdot \  z + y \ \cdot \  z)^\ast  = [t(x \ \cdot \  z )+t(y \ \cdot \  z) - t (x \ \cdot \  z + y \ \cdot \  z)]e + (x \ \cdot \  z)^\ast + (y \ \cdot \  z)^\ast =   [t(x \ \cdot \  z)+t (y \ \cdot \  z) - t(x \ \cdot \  z + y \ \cdot \  z) ]e + z^\ast \ \cdot \  x^\ast  + z^\ast \ \cdot \  y^\ast \ . $
\vskip 2mm
\noindent The beginning and the end of this derivation shows: 
$$\forall (x , y) \in \mathcal{C}^2 \ : \ t(x + y) = t(x) + t(y) \ .$$
The functional $t$ is therefore linear. From (5) then also follows C1:
$$\forall (x , y) \in \mathcal{C}^2 \ : \ (x + y)^\ast =x^\ast + y^\ast \ .$$
Then
$$t(x^\ast) = t(-t(x)e - x) = -t(x) t(e) - t(x) = 2t(x) - t(x)\ ,$$
hence
$$\forall \, x \in \mathcal{C} \ : \ t(x^{\ast}) = t(x) \ ,$$
and from (1) also C4
$$\forall \, x \in \mathcal{C} \ : \ x^{\ast \ast} = x \ .$$

Now let the characteristic of the field\, $\mathbb{F}$\, be\, $\rm{chr}\mathbb{F} \neq 2$ .  Suppose that $z \in \mathcal{A}$ \, is such an element that \, $z = z^\ast$.  Since \ $-t(z)e = z + z^\ast = 2z$\, , \,then \ $z \in \mathbb{F}\{e\}$. We use this knowledge in the following:
$$(x \ \cdot \ x^\ast)^\ast \ = \ x^{\ast \ast} \, \cdot \, x^\ast \ = \ x \, \cdot \, x^\ast \ .$$
Conclusion:
$$\forall \, x \in \mathcal{C} \ : \ x \ \cdot \ x^{\ast} \in \mathbb{F}\{e\} \ $$
and the last axiom C6 is also fulfilled.

By the way, if we define another functional with identity
$$n(x)e := x \ \cdot \ x^\ast \ ,$$
then from (4) we get
$$\forall \, x \in \mathcal{C} \ : \ x \ \cdot \ x + t(x)x + n(x)e \ = \, 0\,. \  \ $$
Algebra with this property is called \textit{quadratic} (sometimes also \textit{conic}).\\

Next suppose that $\rm{chr}\mathbb{F}=2$ . Let \, $(\mathcal{C}, \cdot, e , \ast)$ \,  be an arbitrary commutative algebra in which it holds identically $z^\ast = z$. Then C3 and C5 are fulfilled, but C6 is generally not valid.\\

Let's summarize! The shrunken axiomatic system C3 $\wedge $ C5 is always sufficient, except
\begin{itemize}
    \item  in case of one-dimensional algebra, when all axioms could be replaced solely by C2$^{\rm{var}}$,
    \item  and in case of characteristic 2, when the shrunken system should be C3 $\wedge $  C5 $\wedge $ C6 .
\end{itemize}

\section{}
With all this, however, the task of shrinking axiomatic systems is not yet complete. It has yet to be shown that the proposed contractions are really minimal, so that no proposed axiom is superfluous. This task is not difficult in the four cases above and we shall demonstrate it only in the case of metric.

In fact, here we have three claims:\\
\textbf{M2}$^{\rm{r}}$. \ \ \ \ \ \ \  $d(x , y) \, = \,  0$ \  \ $\Rightarrow \  \ x = y$  \\
\textbf{M2}$^{\rm{1}}$. \ \ \ \ \ \ \  $d(x , y) \, = \,  0$ \  \ $\Leftarrow \  \ x = y$  \\ 
\textbf{M4}$^{\rm{var}}.$ \ \ \ \ \ $d(x , y) \, \leq \, d(x , z) + d(y , z) $ \\
One needs to construct three structures; each should fulfill two of these claims and not the third one. We shall use \ $(\mathcal{M}, d)\, , \ $ where \ $\mathcal{M} = \{a , b\} \, , \ a \neq b$\,, \ and \\
$\neg$ {M2}$^{\rm{r}} \, \wedge$ \, {M2}$^{\rm{1}}$ \, $\wedge \,$ {M4}$^{\rm{var}}$  \ \ \ $d \equiv  0 \ ,$\\
 {M2}$^{\rm{r}} \, \wedge$ \, $\neg$ {M2}$^{\rm{1}}$ \, $\wedge \,$ {M4}$^{\rm{var}}$  \ \ \ $d \equiv   1\ ,$\\
 {M2}$^{\rm{r}} \, \wedge$ \, {M2}$^{\rm{1}}$ \, $\wedge \,$ $\neg$ {M4}$^{\rm{var}}$ \ \ \ $d (a, a) = d(b ,b) = 0 \, , \ d(a , b) = 2d(b , a) = 2   \ .$\\


\begin{thebibliography}{99}

\bibitem{}J. Dieudonné,  \emph{Foundations of Modern Analysis.}  Academic Press, London 1969.
 
\bibitem{} E. O. Kreyszig,     \emph{Advanced Engineering Mathematics} (10th edition). John Wiley \& Sons, New York 2020.

\bibitem{} K. McCrimmon,     \emph{Nonassociative algebras with scalar involution.} Pacific J. Math., Vol. 116, No. 1, 1985.

\bibitem{} J. F. Rigby, J. Wiegold,     \emph{Independent axioms for vector spaces.} Math. Gazette 57.399 (1973), 56-62.

\bibitem{} W. Rudin,    \emph{Real and Complex Analysis.} McGraw-Hill \& Mladinska knjiga, Ljubljana 1970.

\end{thebibliography}
\end{document}